\documentclass[a4paper, 12pt]{article}
\newdimen\paperhight
\usepackage{amssymb,latexsym,amsfonts,amsmath,graphics,oldgerm}

\newcommand{\pf}{\noindent {\it Proof.}\ }
\newcommand{\wt}[1]{\mbox{wt}(#1)}
\newcommand{\mwt}[1]{\mbox{\scriptsize wt}(#1)}

\newtheorem{lemma}{Lemma}
\newtheorem{proposition}[lemma]{Proposition}
\newtheorem{theorem}{Theorem}
\title{Uniform product of 
$A_{g,n}(V)$ for an orbifold model 
$V$ and $G$-twisted Zhu algebra}
\author{
\begin{tabular}{c}
Masahiko Miyamoto
\footnote{Supported by the Grants-in-Aids for 
Scientific Research, 
No. 13440002 and No. 12874001, The Ministry of 
Education, Science and 
Culture, Japan.}
\\
and 
\\
Kenichiro Tanabe\\
\\
Institute of Mathematics \\
University of Tsukuba
\end{tabular}}
\date{February 12, 2002}
\begin{document}
\maketitle
\begin{abstract}
Let $V$ be a vertex operator algebra and $G$ a finite 
automorphism group of $V$. For each $g\in G$ and 
nonnegative rational number 
$n\in {\mathbb Z}/|g|$, a $g$-twisted Zhu algebra 
$A_{g,n}(V)$ plays an important role in the theory
of vertex operator algebras, 
but the given product in 
$A_{g,n}(V)$ depends on the eigenspaces of $g$. 
We show that there is a uniform definition of 
products on $V$ and we introduce a $G$-twisted 
Zhu algebra $A_{G,n}(V)$
which covers all $g$-twisted Zhu algebras. 

Let $V$ be simple and let ${\cal S}$ 
be a finite set of inequivalent irreducible twisted 
$V$-modules which is closed under the action of $G$. 
There is a finite dimensional
semisimple associative algebra 
${\cal A}_{\alpha}(G,{\cal S})$ for
a suitable 2-cocycle naturally determined by the 
$G$-action on ${\cal S}$.
We show that a duality theorem of Schur-Weyl type
holds for the actions of ${\cal A}_{\alpha}(G,{\cal S})$ and 
$V^G$ on the direct sum of twisted $V$-modules in ${\cal S}$
as an application of the theory of $A_{G,n}(V)$.
It follows as a natural consequence of the result
that for any $g\in G$ every irreducible $g$-twisted 
$V$-module is a completely reducible $V^G$-module.
\end{abstract}

\section{Introduction}
\setcounter{equation}{0}
Let $V$ be a vertex operator algebra 
(cf. \cite{B},\cite{FHL},\cite{FML})
and $G$ a finite automorphism group of $V$ order $T$. 
For an investigation of $V$-modules, 
a Zhu algebra $A(V)=V/O(V)$ (or its extension 
$A_n(V)=V/O_n(V)$),
which was introduced in \cite{Z} and \cite{DLM1-5},
plays an important role. For a $g$-twisted $V$-module, 
a $g$-twisted Zhu algebra 
$A_{g,n}(V)=V/O_{g,n}(V)$ was introduced in \cite{DLM2} and 
plays a similar role, where nonnegative 
$n\in {\mathbb Z}/T$. 
However, the definition of product $u\ast_{g,n} v$ 
in $A_{g,n}(V)$ 
depends on the choice of eigenspaces of $g$ containing $u$.  
Namely, for $g\in G$, $V$ decomposes into the direct sum of 
eigenspaces $V=\oplus_{r=0}^{T-1} V^{(g,r)}$, 
where $V^{(g,r)}$ is the eigenspace of $g$ with eigenvalue 
$e^{-2\pi \sqrt{-1}r/T}$. 
If $u\in V^{(g,0)}$, then 
\[u\ast_{g,n}v\equiv \sum_{m=0}^{l}(-1)^m{m+h\choose l}
\mbox{Res}_zY(u,z)v\frac{(1+z)^{\mwt{u}+l}}{z^{l+m+1}} 
\pmod{O_{g,n}(V)}\]
and if $u\in V^{(g,r)}$ $(r\not=0)$, 
then $u\ast_{g,n}v\equiv 0 \pmod{O_{g,n}(V)}$, 
where $n=l+i/T$ with $0\leq i\leq T-1$ and 
$l\in {\mathbb Z}_{\geq 0}$. 
One of the purposes in this paper is to show that for 
$u,v\in V$, 
there is a unified form $\ast_n$ of products which does 
only depend on the weights of $u$ and $v$ such that 
\[u\ast_nv\equiv u\ast_{g,n} v \ \pmod{O_{g,n}(V)}\]
for any $g\in G$. Using this product, a natural map 
\[\phi_n: V \to \oplus_{g\in G}A_{g,n}(V)\]
given by $\phi_n(u)=(u,u,\cdots,u)$ becomes an 
algebra homomorphism:
\[\phi_n(u\ast_n v)=(u \ast_{g,n}v+O_{g,n}(V))_{g\in G}.\]  
Set $A_{G,n}(V)=V/\cap_{g\in G}O_{g,n}(V)$.
We will call $(A_{G,n}(V),*_n)$ a $G$-twisted Zhu algebra. 
We will show that if $V$ is simple, 
then $A_{G,n}(V)\simeq \oplus_{g\in G}A_{g,n}(V)$ as algebras.

We will apply the theory of $A_{G,n}(V)$ to 
a problem of the representation theory of $V^G$.
One of the important problems in the orbifold conformal field 
theory is to determine the $V^G$-modules. 
In particular, it is one of the main conjectures that 
if $V$ is rational, then $V^G$ is rational, that is, 
every $V^G$-module 
is completely reducible (cf. \cite{DVVV}). 
For a simple VOA $V$, it is shown in \cite{DY} 
that every irreducible $V$-module is 
a completely reducible $V^G$-module as a natural consequence
of a duality theorem of Schur-Weyl type.
Another important category is a twisted module 
(cf. \cite{DLM1}). 
Namely, for any $g\in G$, any $g$-twisted module is also 
a $V^G$-module. 
The same results as those in the untwisted case 
are shown in \cite{Y} for irreducible
$g$-twisted $V$-modules when $g$ is in the center of $G$. 

In this paper we extend their results to any $g\in G$
and any irreducible $g$-twisted $V$-module.
Let's state our results more explicitly. 
There is a natural
right $G$-action on the set of all inequivalent 
irreducible twisted $V$-modules.
Let ${\cal S}$ 
be a finite set of inequivalent irreducible twisted
$V$-modules which is closed under the action of $G$. 
We define a 
finite dimensional semisimple associative algebra 
${\cal A}_{\alpha}(G,{\cal S})$ over ${\mathbb C}$
associated to $G$, ${\cal S}$ and a suitable $2$-cocycle
$\alpha$ naturally determined by the $G$-action 
on ${\cal S}$. The algebra ${\cal A}_{\alpha}(G,{\cal S})$
is constructed in more general setting in \cite{DY}
and is called a generalized twisted double.
We show a duality theorem of Schur-Weyl type
for the actions 
of $V^G$ and ${\cal A}_{\alpha}(G,{\cal S})$ on 
the direct sum of twisted $V$-modules
in ${\cal S}$ which is denoted by ${\cal M}$.
That is, each simple ${\cal A}_{\alpha}(G,{\cal S})$ 
occurs in ${\cal M}$
and its multiplicity space is an irreducible $V^G$-module. 
Moreover, the different multiplicity spaces are inequivalent 
$V^G$-modules.
It follows as a natural consequence of the result
that for any $g\in G$ every irreducible $g$-twisted 
$V$-module is a completely reducible $V^G$-module.
The theory of $A_{G,n}(V)$ 
allows us to reduce an infinite dimensional problem
to a finite dimensional one.

This paper is organized as follows. In Section 2, we first
review the Zhu algebras of a vertex operator algebra. 
We introduce $(A_{G,n}(V), *_n)$ and study their properties.
In Section 3, we define a
finite dimensional semisimple associative algebra 
${\cal A}_{\alpha}(G,{\cal S})$ over ${\mathbb C}$ and
show a duality theorem of Schur-Weyl type
as an application of the theory of $A_{G,n}(V)$.
In Section 4, we compute the determinant of a
matrix used in Section 2.

\section{Existence of unified form of product and 
\setcounter{equation}{0}
$G$-twisted Zhu algebras}
We fix some notation which will be in force throughout 
the paper. $V$ is a vertex operator 
algebra and $G$ is a finite automorphism group of 
$V$ of order $T$. 
For $g\in G$, set $V^{(g,r)}
=\{u\in V|\ gu=e^{-2\pi\sqrt{-1}r/T}u\}$ for $0\leq r\leq T-1$.
For $u\in V$ we denote by $u^{(g,r)}$ the $r$-th 
component of $u$ in the decomposition 
$V=\oplus_{r=0}^{T-1}V^{(g,r)}$, that is, 
$u=\sum_{r=0}^{T-1}u^{(g,r)}$,
$u^{(g,r)}\in V^{(g,r)}\ (0\leq r\leq T-1)$.

In this section, our main purpose is to 
show that there is a unified product form 
$u\ast_n v$ which depends only on the weights of $u$ and $v$.
We introduce an associative algebra $(A_{G,n}(V),*_n)$ of $V$ 
associated to $G$ and nonnegative $n\in {\mathbb Z}/T$
and study their properties.

We first recall the Zhu algebra $A_{g,n}(V)$ of $V$
associated with $g\in G$ and nonnegative $n\in {\mathbb Z}/T$
introduced in \cite{DLM2}. This
algebra was first introduced in \cite{Z}
for the case when $g$ is the identity element and $n=0$.
Fix nonnegative 
$n=l+i/T\in {\mathbb Z}/T$ with $l$ a nonnegative integer and
$0\leq i\leq T-1$.

For $0\leq r\leq T-1$ we define 
$\delta_i(r)=1$ if $i\geq r$ and $\delta_i(r)=0$ if $i<r$.
We also set $\delta_i(T)=1$.
Let $g\in G$.
Let $O_{g,n}(V)$ be the linear span of all $u\circ_{g,n}v$ and 
$(L(-1)+L(0))v$ where for homogeneous $u\in V$ and $v\in V$,
$u\circ_{g,n}v$ denotes 
\[\sum_{r=0}^{T-1}\mbox{Res}_zY(u^{(g,r)},z)v
\frac{(1+z)^{\mwt{u}-1+\delta_i(r)+l+r/T}}
{z^{2l+\delta_i(r)+\delta_i(T-r)}}.\]
Define the linear space $A_{g,n}(V)$ to be the quotient 
$V/O_{g,n}(V)$.
We also define a second product $*_{g,n}$ on $V$ for 
homogeneous 
$u\in V$ and $v\in V$ by 
\begin{eqnarray*}
u*_{g,n}v & = & \sum_{m=0}^{l}(-1)^m{m+l\choose l}
\mbox{Res}_zY(u^{(g,0)},z)v
\frac{(1+z)^{\mwt{u}+l}}{z^{l+m+1}}\\
& = & 
\sum_{k=0}^{\infty}
\sum_{j=0}^{k}(-1)^{l-j}{2l-j\choose l}
{\wt{u}+l\choose k-j}u^{(g,0)}_{-2l-1+k}v.
\end{eqnarray*}
Extend this linearly to obtain a bilinear product on $V$.
Note that if $i=j+kT/|g|$ where $0\leq j\leq T/|g|-1$ and
$0\leq k\leq |g|-1$, then $A_{g,n}(V)=A_{g,kT/|g|}(V)$. 
We recall the following properties of $O_{g,n}(V)$ 
and $A_{g,n}(V)$.

\begin{lemma}\label{Lemma : ele}
{\rm (\cite{DLM2}, Lemma 2.1, Lemma 2.2, Theorem 2.4, 
and Proposition 2.5)}
\begin{enumerate}
\item $\oplus_{r=1}^{T-1}V^{(g,r)}\subset O_{g,n}(V)$.
\item For homogeneous $u\in V$, $v\in V$ and integers 
$0\leq k\leq m$,
\[\sum_{r=0}^{T-1}\mbox{\rm Res}_zY(u^{(g,r)},z)v
\frac{(1+z)^{\mwt{u}-1+\delta_i(r)+l+r/T+k}}
{z^{2l+\delta_i(r)+\delta_i(T-r)+m}}\in O_{g,n}(V).\]
\item $(A_{g,n}(V), *_{g,n})$ is an associative algebra 
with the identity element. 
\item The identity map on $V$ induces an onto algebra 
homomorphism from $A_{g,n}(V)$ to $A_{g,n-1/T}(V)$. 
\end{enumerate}
\end{lemma}
Set $A_{G,n}(V)=V/\cap_{g\in G}O_{g,n}(V)$.
We state our main result in this section.

\begin{theorem}\label{Theorem : zhu}
\begin{enumerate}
\item For $u,v\in V$,
there exists a unified product $u*_n v$
for $A_{g,n}(V)$ for any $g\in G$.
That is,
\[u*_nv\equiv u*_{g,n}v \pmod{O_{g,n}(V)}\]
for any $g\in G$. In particular, $(A_{G,n}(V),*_n)$
is an associative algebra.
Moreover, if $V$ is simple, then the injective map 
$\varphi_n : A_{G,n}(V)\rightarrow \oplus_{g\in G}A_{g,n}(V)$
defined by 
$\varphi_n(u+\cap_{g\in G}O_{g,n}(V))=\sum_{g\in G}
(u+O_{g,n}(V))$ for $u\in V$ is an onto algebra homomorphism.
That is, $A_{G,n}(V)\simeq \oplus_{g\in G}A_{g,n}(V)$ 
as algebras. 
\item
The identity map on $V$ induces 
an onto algebra homomorphism
{from} $A_{G,n}(V)$ to $A_{G,n-1/T}(V)$. 
\end{enumerate}
\end{theorem}

\noindent{\bf Remark.}
We will construct $u*_nv$ as a linear combination of 
$\{u_jv\ |\ j\in {\mathbb Z}\}$
in the proof of Theorem \ref{Theorem : zhu}:
$u*_nv=\sum_{j\in {\mathbb Z}}\gamma_ju_jv$
where 
$\gamma_{j}\in {\mathbb C}\ (j\in {\mathbb Z})$ 
are all zero except finite numbers.
Then the set $\{\gamma_j\}_{j\in {\mathbb Z}}$ 
depends only on $n$ and the weights of $u$ and $v$.

\medskip
We will call $A_{G,n}(V)$ a $G$-twisted Zhu algebra.
We need some lemmas in order to show Theorem 
\ref{Theorem : zhu}.
Let $u,v\in V$ and fix $M\in {\mathbb Z}$ with
$M\geq \wt{u}+\wt{v}-1$ in the next three lemmas. 
Note that 
$u^{(g,r)}_jv=0$ for all $g\in G, 0\leq r\leq T-1$, 
and $j>M$.
Set $Q_r=\wt{u}-1+\delta_i(r)+l+r/T$ for $0\leq r\leq T-1$.

\begin{lemma}\label{Lemma : mod}
For $q\in {\mathbb Z}$ with $-2l-1\geq q$ and $m=1,2,\ldots$,
we have
\begin{eqnarray*}
u_{q-m}v & \equiv & 
\sum_{k=0}^{Q_0-1}\sum_{j=0}^{k}{Q_0\choose j}
{-Q_0\choose k+m-j}u^{(g,0)}_{q+k}v\\
& & +\sum_{r=1}^{T-1}
\sum_{k=0}^{M-q}\sum_{j=0}^{k}{Q_r\choose j}
{-Q_r\choose k+m-j}u^{(g,r)}_{q+k}v,\pmod{O_{g,n}(V)},
\end{eqnarray*}
for any $g\in G$.
\end{lemma}
\pf 
Let $g\in G$. 
We may assume that $u$ is homogeneous and $u\in V^{(g,r)}$ 
for $r\in \{0,\ldots,T-1\}$. 
Since $2\geq \delta_i(r)-\delta_i(T-r)$, we have
\begin{eqnarray*}
O_{g,n}(V) & \ni & 
\mbox{Res}_zY(u,z)v\frac{(1+z)^{Q_r}}{z^{2l+2+m}}\\
& = & \sum_{j=0}^{\infty}{Q_r\choose j}u_{-2l-2-m+j}v.  
\end{eqnarray*}
for any nonnegative integer $m$ {from} 
Lemma \ref{Lemma : ele} (2).
Let $q\in {\mathbb Z}$ with $q\leq -2l-1$.
Since $q\leq -2l-1$, we have 
\begin{eqnarray*}
O_{g,n}(V) & \ni & 
\sum_{j=0}^{\infty}{Q_r\choose j}u_{q-m+j}v \\
& = & 
\sum_{j=0}^{m-1}{Q_r\choose j}u_{q-m+j}v+
\sum_{k=0}^{\infty}{Q_r\choose m+k}u_{q+k}v
\end{eqnarray*}
for any positive integer $m$. Namely, we have
\begin{eqnarray}\label{eqn : sankaku}
\begin{array}{rcl}
u_{q-1}v & \equiv & -{Q_r\choose 1}u_{q}v
-{Q_r\choose 2}u_{q+1}v-\cdots,\\
u_{q-2}v+{Q_r\choose 1}u_{q-1}v & \equiv & 
-{Q_r\choose 2}u_{q}v-{Q_r\choose 3}u_{q+1}v-\cdots,\\
u_{q-3}v+{Q_r\choose 1}u_{q-2}v+{Q_r\choose 2}u_{q-1}v 
& \equiv & 
-{Q_r\choose 3}u_{q}v-{Q_r\choose 4}u_{q+1}v-\cdots,\\
\multicolumn{1}{c}{\vdots} & & \multicolumn{1}{c}{\vdots}
\end{array}
\end{eqnarray}
module $O_{g,n}(V)$. 
Solving these congruent equations, 
$u_{q-m}v$ is congruent to a linear combination of 
$\{ u_{q}v, u_{q+1}v,\ldots,u_{M}v\}$ modulo $O_{g,n}(V)$ 
for any positive integer $m$.     
If $r=0$, 
then $Q_r=\wt{u}+l \in {\mathbb Z}_{\geq 0}$ and so 
$u_{q-m}v$ is congruent to 
a linear combination of 
$\{u_{q}v, u_{q+1}v,\ldots,u_{q+Q_0-1}v\}$.

Replacing $u_jv$ by $x^j$,  
we will investigate the coefficients by viewing 
the above equations as a Lorentz series. 
Set 
\begin{eqnarray}
f_r(x) & = & \sum_{j=0}^{\infty}{Q_r\choose j}x^j.
\label{eqn : formal}
\end{eqnarray}
We note that if $\sum_{j\in {\mathbb Z}}
\alpha_jx^j\in 
{\mathbb C}[x^{-1}]x^{-2l-2}f_r(x)$, then 
$\sum_{j\in {\mathbb Z}} \alpha_ju_jv\in O_{g,n}(V)$ 
and 
$f_r(x)$ depends only on $N=\wt{u}$, $r$ and $n=l+i/T$.   
Set $O^r(N:x)=
{\mathbb C}[x^{-1}]x^{-2l-2}f_r(x)\subset 
{\mathbb C}[x^{-1}][[x]]$. 

Let $m$ be a positive integer. In order to express
$u_{q-m}v$ by a linear combination of 
$\{u_{q}v, u_{q+1}v,\ldots, u_{M}v\}$ 
modulo $O_{g,n}(V)$, we first write $1/f_r(x)$ as a sum 
\[\frac{1}{f_r(x)}=(f_r(x)^{-1})_{<m}
+(f_r(x)^{-1})_{\geq m}, \] 
where $(f_r(x)^{-1})_{<m}$ is 
the part whose terms have degree less than $m$ 
and $(f_r(x)^{-1})_{\geq m}$ is the part whose terms have 
degree greater than or equal to $m$. 
In particular, we have 
$(f_r(x)^{-1})_{\geq m}\in x^{m}{\mathbb C}[[x]]$. 
Since $(x^{q-1}f_r(x))(f_r(x)^{-1})_{<m}/x^{m-1}\equiv 0 
\pmod{O^r(N:x)}$, we obtain 
\begin{eqnarray*}
x^{q-m} & \equiv& x^{q-m}-x^{q-1}f_r(x)
\frac{(f_r(x)^{-1})_{<m}}{x^{m-1}}\\
& \equiv & x^{q-m}-x^{q-1}f_r(x)
\frac{f_r(x)^{-1}-(f_r(x)^{-1})_{\geq m}}{x^{m-1}}\\
& \equiv & x^{q-m}f_r(x)(f_r(x)^{-1})_{\geq m}\\
& \equiv & x^{q-m}(\sum_{j=0}^{\infty}{Q_r\choose j}x^j)
(\sum_{k=m}^{\infty}{-Q_r\choose k}x^k)\\
& \equiv & 
\sum_{k=0}^{\infty}
\sum_{j=0}^{k}{Q_r\choose j}{-Q_r\choose k+m-j}
x^{q+k}\pmod{O^r(N:x)}. 
\end{eqnarray*}
Hence we have
\begin{eqnarray*}
u_{q-m}v & \equiv & \sum_{k=0}^{\infty}
\sum_{j=0}^{k}{Q_r\choose j}{-Q_r\choose k+m-j}
u_{q+k}v\\
& \equiv & \sum_{k=0}^{M-q}
\sum_{j=0}^{k}{Q_r\choose j}{-Q_r\choose k+m-j}
u_{q+k}v \pmod{O_{g,n}(V)}. 
\end{eqnarray*}
$\Box$

\begin{lemma}\label{Lemma : kiri}
For $r\in \{0,\ldots,T-1\}$ and $p\in {\mathbb Z}$,
there exists $F^{r,p}(u,v,n)\in V$ such that
\begin{eqnarray}
F^{r,p}(u,v,n)\equiv u^{(g,r)}_pv\ \pmod{O_{g,n}(V)}
\label{eqn : mod}
\end{eqnarray}
for any $g\in G$.
\end{lemma}
\pf We may assume $p\leq M$.
Define a set 
\[S=\left\{(s,k)\left|
\begin{array}{ll}
s=0,1,\ldots,T-1\\
k=0,\ldots,Q_0-1 & \mbox{ for }s=0 \qquad \mbox{  and }\cr
k=0,\ldots,M-q &\mbox{ for }s=1,\ldots,T-1   
\end{array}\right. \right\}.
\]
Let $q=\mbox{min}\{p,-2l-1\}$ and set
\begin{eqnarray*}
F^{r,p}(u,v,n) & = & 
\left\{
\begin{array}{ll}
{\displaystyle
u_pv-\sum_{m=1}^{|S|}
\lambda_{q-m}u_{q-m}v} & \mbox{if }r=0,\\
{\displaystyle
\sum_{m=1}^{|S|}\lambda_{q-m}u_{q-m}v} 
& \mbox{if }1\leq r\leq T-1,
\end{array}
\right.
\end{eqnarray*}
where $\lambda_{q-m}\in {\mathbb C}\ (m=1,2,\ldots,|S|)$.
We show that it is possible
to choose a set $\{\lambda_{q-m}\}_{m=1}^{|S|}$
which depends only on $\wt{u},\wt{v}$ and $n$
such that the the assertion holds.
For $g\in G$, we have
\begin{eqnarray*}
& & \sum_{m=1}^{|S|}\lambda_{q-m}u_{q-m}v\\
& \equiv & \sum_{m=1}^{|S|}\lambda_{q-m}
\Big(\sum_{k=0}^{Q_0-1}\sum_{j=0}^{k}{Q_0\choose j}
{-Q_0\choose k+m-j}u^{(g,0)}_{q+k}v\\
& & +\sum_{s=1}^{T-1}\sum_{k=0}^{M-q}\sum_{j=0}^{k}
{Q_s\choose j}
{-Q_s\choose k+m-j}u^{(g,s)}_{q+k}v\Big)\\
& \equiv & \sum_{k=0}^{Q_0-1}
\sum_{m=1}^{|S|}\lambda_{q-m}\sum_{j=0}^{k}
{Q_0\choose j}{-Q_0\choose k+m-j}u^{(g,0)}_{q+k}v\\
& & +\sum_{s=1}^{T-1}
\sum_{k=0}^{M-q}\sum_{m=1}^{|S|}\lambda_{q-m}
\sum_{j=0}^{k}{Q_s\choose j}
{-Q_s\choose k+m-j}u^{(g,s)}_{q+k}v\ \pmod{O_{g,n}(V)}
\end{eqnarray*}
{from} Lemma \ref{Lemma : mod}.
Comparing both sides of formula (\ref{eqn : mod}),
we have  $|S|$ linear equations:
\begin{equation}
\label{eqn : linear}
\begin{array}{rcl}
\multicolumn{3}{l}{{\bullet} 
\mbox{ In the case $1\leq r\leq T-1$,}}\\
\sum_{m=1}^{|S|}\lambda_{q-m}
\sum_{j=0}^{k}{Q_s\choose j}
{-Q_s\choose k+m-j} 
& = & 
\left\{
\begin{array}{ll}
1 & \mbox{if $s=r$ and $q+k=p$},\\
0 & \mbox{otherwise.}
\end{array}\right.\\
\multicolumn{3}{l}{\bullet \mbox{ In the case $r=0$,}}\\
\sum_{m=1}^{|S|}\lambda_{q-m}
\sum_{j=0}^{k}{Q_s\choose j}
{-Q_s\choose k+m-j} 
& = & 
\left\{
\begin{array}{ll}
1 & \mbox{if $s\neq 0$ and $q+k=p$},\\
0 & \mbox{otherwise.}
\end{array}\right.
\end{array}
\end{equation}
Here $(s,k)$ ranges over $S$.
We denote $\sum_{j=0}^{k}{Q_s\choose j}{-Q_s\choose k+m-j}$ 
by $\alpha^{s,k}_m$ for $(s,k)\in S$ and a 
positive integer $m$.
We set a $|S|\times |S|$-matrix
\begin{eqnarray*}
A_1 & = & \left[
\begin{array}{ccccccc}
\alpha^{0,0}_{1}&\cdots&\alpha^{0,Q_0-1}_1&
\alpha^{1,0}_{1}&\cdots &\cdots 
&\alpha^{T-1,M-p}_{1}\\
\alpha^{0,0}_{2}&\cdots&\alpha^{0,Q_0-1}_2&
\alpha^{1,0}_{2}&\cdots &\cdots 
&\alpha^{T-1,M-p}_{1}\\
\vdots & &\vdots&\vdots& & &\vdots \\
\alpha^{0,0}_{|S|}&\cdots&\alpha^{0,Q_0-1}_{|S|}&
\alpha^{1,0}_{|S|}&\cdots &\cdots 
&\alpha^{T-1,M-p}_{|S|}\\
\end{array}\right]. 
\end{eqnarray*}
It is sufficient to show that the matrix $A_1$ is non-singular
in order to prove the equations (\ref{eqn : mod}) have
a solution $\{\lambda_{q-m}\}_{m=1}^{|S|}$.

We set a $|S|\times |S|$-matrix
\begin{eqnarray*}
A_2 & = & \left[
\begin{array}{ccccccc}
{-Q_0\choose 1}&\cdots&{-Q_0\choose Q_0}&
{-Q_1\choose 1}&\cdots &\cdots & 
{-Q_{T-1}\choose M-q+1}\\
{-Q_0\choose 2}&\cdots&{-Q_0\choose Q_0+1}&
{-Q_1\choose 2}&\cdots &\cdots & 
{-Q_{T-1}\choose M-q+2}\\
\vdots & &\vdots&\vdots& & &\vdots \\
{-Q_0\choose |S|}&\cdots&{-Q_0\choose Q_0+|S|-1}&
{-Q_1\choose |S|}&\cdots &\cdots & 
{-Q_{T-1}\choose M-q+|S|}
\end{array}\right]. 
\end{eqnarray*}
Then $\det A_1=\det A_2$ because
\begin{eqnarray*}
& & 
(\alpha^{s,0}_m,\alpha^{s,1}_m\ldots,\alpha^{s,M-q}_m)\\
& = & ({-Q_s\choose m},{-Q_s\choose m+1},\ldots,
{-Q_s\choose m+M-q})
\left[
\begin{array}{ccccc}
1 & {Q_s\choose 1} & {Q_s\choose 2} & \cdots & 
{Q_s\choose M-q}\\
0 & 1 & {Q_s\choose 1} & \ddots & \vdots\\
\vdots & \ddots & \ddots & \ddots & {Q_s\choose 2}\\
\vdots & & \ddots & \ddots & {Q_s\choose 1}\\
0 & \cdots & \cdots & 0 & 1
\end{array}
\right]
\end{eqnarray*}
for $1\leq s\leq T-1$ and a similar formula for $s=0$.
So it sufficients to show that the matrix $A_2$ is 
non-singular.
The following computation is used in 
a transformation of $A_2$. For $(s,k)\in S$, we obtain 
\begin{eqnarray*}
& & \sum_{m=1}^{\infty}
\sum_{j=0}^{m-1}{Q_0\choose m-1-j}{-Q_s\choose k+1+j}x^{m}\\
& = & \sum_{m=0}^{\infty}
\sum_{j=0}^{m}{Q_0\choose m-j}{-Q_s\choose k+1+j}x^{m+1}\\
& = & \sum_{m=0}^{\infty}
\sum_{j=0}^{m}{Q_0\choose j}{-Q_s\choose k+m+1-j}x^{m+1}\\
& = & x^{-k}
(1+x)^{Q_0}((1+x)^{-Q_s}
-\sum_{j=0}^{k}{-Q_s\choose j}x^j)\\
& = & x^{-k}
(1+x)^{Q_0-Q_s}-\sum_{t=0}^{Q_0}\sum_{j=0}^{k}
{Q_0\choose t}{-Q_s\choose j}x^{t+j-k}\\
& = & \sum_{m=0}^{\infty}
{Q_0-Q_s\choose m}x^{m-k}
-\sum_{m=0}^{Q_0+k}\sum_{j=0}^{k}
{Q_0\choose m-j}{-Q_s\choose j}x^{m-k}\\
& = & \sum_{m=1}^{\infty}
{Q_0-Q_s\choose m+k}x^{m}
-\sum_{m=1}^{Q_0}\sum_{j=0}^{k}
{Q_0\choose k+m-j}{-Q_s\choose j}x^{m}.
\end{eqnarray*}
So we have
\begin{eqnarray}
& & \left[
\begin{array}{cccc}
1 & 0 & \cdots & 0\\
{Q_0\choose 1} & \ddots & \ddots & \vdots\\
\vdots & \ddots & \ddots & 0\\
{Q_0\choose |S|-1} & \cdots & {Q_0\choose 1} & 1\\
\end{array}
\right]
\left[
\begin{array}{cccc}
{-Q_0\choose 1} & {-Q_0\choose 2} & 
\cdots & {-Q_0\choose Q_0}\\
{-Q_0\choose 2} & {-Q_0\choose 3} & 
\cdots & {-Q_0\choose Q_0+1}\\
\vdots & \vdots & & \vdots \\
{-Q_0\choose |S|} & {-Q_0\choose |S|+1} & \cdots & 
{-Q_0\choose Q_0+|S|-1}
\end{array}\right] \nonumber\\
& = & 
-\left[
\begin{array}{ccccc}
{Q_0\choose 1} & \cdots & {Q_0\choose Q_0-1} & 1\\
\vdots & \mbox{\reflectbox{$\ddots$}} 
& \mbox{\reflectbox{$\ddots$}} & 0 \\
{Q_0\choose Q_0-1} & \mbox{\reflectbox{$\ddots$}} & 
\mbox{\reflectbox{$\ddots$}} & \vdots \\
1 & 0 & \cdots & 0\\\hline
0 &\cdots  & \cdots & 0\\
\vdots & & & \vdots\\
0 & \cdots & \cdots & 0\\
\end{array}\right]
\left[
\begin{array}{cccc}
1 & {Q_0\choose 1} & \cdots & {Q_0\choose Q_0-1}\\
0 & \ddots & \ddots & \vdots\\
\vdots & \ddots & \ddots & {Q_0\choose 1}\\
0 & \cdots & 0 & 1
\end{array}
\right]\label{eqn : det0}
\end{eqnarray}
and for $1\leq s\leq T-1$,
\begin{eqnarray}
& & \left[
\begin{array}{cccc}
1 & 0 & \cdots & 0\\
{Q_0\choose 1} & \ddots & \ddots & \vdots\\
\vdots & \ddots & \ddots & 0\\
{Q_0\choose |S|-1} & \cdots & {Q_0\choose 1} & 1\\
\end{array}
\right]
\left[
\begin{array}{cccc}
{-Q_s\choose 1} & {-Q_s\choose 2} & 
\cdots & {-Q_s\choose M-p+1}\\
{-Q_s\choose 2} & {-Q_s\choose 3} & 
\cdots & {-Q_s\choose M-p+2}\\
\vdots & \vdots & & \vdots \\
{-Q_s\choose |S|} & {-Q_s\choose |S|+1} & \cdots & 
{-Q_s\choose M-p+|S|}
\end{array}\right] \nonumber\\
& = & 
\left[
\begin{array}{ccccc}
y_{10} & y_{11} & \cdots & y_{1,M-p}\\
\vdots & \vdots & & \vdots \\
y_{Q_0,0} & y_{Q_0,1} & \cdots & y_{Q_0,M-p}\\\hline
{Q_0-Q_s\choose Q_0+1} & {Q_0-Q_s\choose Q_0+2} & \cdots & 
{Q_0-Q_s\choose Q_0+M-p+1}\\
{Q_0-Q_s\choose Q_0+2} & {Q_0-Q_s\choose Q_0+3} & \cdots & 
{Q_0-Q_s\choose Q_0+M-p+2}\\
\vdots & \vdots & & \vdots\\
{Q_0-Q_s\choose |S|} & {Q_0-Q_s\choose |S|+1} & \cdots & 
{Q_0-Q_s\choose M-p+|S|}\\
\end{array}\right]
\label{eqn : dets}
\end{eqnarray}
where $y_{m,k}$ denotes ${Q_0-Q_s\choose m+k}-
\sum_{j=0}^{k}{Q_0\choose m+k-j}{-Q_s\choose j}$
for $1\leq m\leq Q_0$ and $0\leq k\leq M-p$.
We hence have $\det A_2=(-1)^{Q_0(Q_0+1)/2}\det A_3$, where 
$A_3$ is the following $(|S|-Q_0)\times (|S|-Q_0)$-matrix:
\begin{eqnarray}
& & A_3= \nonumber\\
&  & \left[
\begin{array}{cccccc}
{Q_0-Q_1\choose Q_0+1} & \cdots & {Q_0-Q_1\choose Q_0+M-p+1}
&{Q_0-Q_2\choose Q_0+1} & \cdots & 
{Q_0-Q_{T-1}\choose Q_0+M-p+1}\\
{Q_0-Q_1\choose Q_0+2} & \cdots & {Q_0-Q_1\choose Q_0+M-p+2}&
{Q_0-Q_2\choose Q_0+2} & \cdots & 
{Q_0-Q_{T-1}\choose Q_0+M-p+2}\\
\vdots & & \vdots & \vdots & & \vdots \\
{Q_0-Q_1\choose |S|} & \cdots & {Q_0-Q_1\choose M-p+|S|}
& {Q_0-Q_2\choose |S|} & \cdots & 
{Q_0-Q_{T-1}\choose M-p+|S|}
\end{array}
\right].
\label{eqn : det2}
\end{eqnarray}
It is proved in Appendix that $A_3$ is non-singular.
$\Box$

\begin{lemma}\label{Lemma : wa}
Let $V$ be simple and let $g,h\in G$ with $g\neq h$.
Then $O_{g,n}(V)+O_{h,n}(V)=V$.
\end{lemma}
\pf
For $g,h\in G$ with $g\neq h$, 
there are $k,r \in \{0,\ldots,T-1\}$ with $k\neq r$
and $v\in V$ such that 
$0\neq v\in V^{(g,k)}$ and $v^{(h,r)}\neq 0$.
For any $u\in V$ and $p\in {\mathbb Z}$, 
consider $F^{-r,p}(u,v,n)$ in Lemma \ref{Lemma : kiri}.
Then 
\[F^{-r,p}(u,v,n)\equiv u^{(g,-r)}_pv\equiv 0
\pmod{O_{g,n}(V)}\]
{from} Lemma \ref{Lemma : ele} (1)
because $u^{(g,-r)}_pv\in V^{(g,k-r)}$ 
and $k-r\not\equiv 0\pmod{T}$.
We also have 
\begin{eqnarray*}
F^{-r,p}(u,v,n) & \equiv & u^{(h,-r)}_pv\\
& \equiv & u^{(h,-r)}_pv^{(h,r)}\\
& \equiv & u_pv^{(h,r)}\pmod{O_{h,n}(V)}.\
\end{eqnarray*}
We hence have $u_pv^{(h,r)}\in O_{g,n}(V)+O_{h,n}(V)$ 
for any $u\in V$ and $p\in {\mathbb Z}$. 
Since $v^{(h,r)}\neq 0$ and $V$ is simple,
$V=\mbox{Span}\{u_pv^{(h,r)}\ |\ u\in V, p\in {\mathbb Z}\}$
{from} Proposition 4.1 in \cite{DM}.
So $O_{g,n}(V)+O_{h,n}(V)=V$ holds.
$\Box$

\medskip
Now we start to prove Theorem \ref{Theorem : zhu}.
Let $u,v\in V$ and fix $M\in {\mathbb Z}$ with
$M\geq \wt{u}+\wt{v}-1$.
Define 
\begin{eqnarray*}
u*_nv & = & \sum_{k=0}^{M+2l+1}
\sum_{j=0}^{k}(-1)^{l-j}{2l-j\choose l}
{\wt{u}+l\choose k-j}F^{-2l-1+k,0}(u,v,n).
\end{eqnarray*}
We have $u*_nv\equiv u*_{g,n}v\pmod{O_{g,n}(V)}$ 
for any $g\in G$
{from} Lemma \ref{Lemma : kiri}.
If $V$ is simple, then we have 
$A_{G,n}(V)\simeq \oplus_{g\in G}A_{g,n}(V)$ as algebras 
using the Chinese remainder theorem and Lemma \ref{Lemma : wa}.
So (1) holds. (2) is clear {from} (1) and Lemma 
\ref{Lemma : ele}.
$\Box$

\section{A duality theorem of Schur-Weyl type}
\setcounter{equation}{0}
We always assume that $V$ is simple throughout this section.
A finite dimensional 
semisimple associative algebra 
${\cal A}_{\alpha}(G, {\cal S})$ over ${\mathbb C}$
associated to $G$,
a finite right $G$-set ${\cal S}$ 
and a suitable $2$-cocycle $\alpha$ is 
constructed in \cite{DY}.
${\cal A}_{\alpha}(G, {\cal S})$ is called the generalized
twisted double there.
In this section we first review its
construction in 
the case that ${\cal S}$ is a finite $G$-stable set 
of inequivalent irreducible twisted $V$-modules
and a $2$-cocycle $\alpha$ 
naturally determined by the $G$-action on ${\cal S}$.
We will show a duality theorem of 
Schur-Weyl type for the actions of $V^G$ and 
${\cal A}_{\alpha}(G, {\cal S})$ on 
${\cal M}=\oplus_{(g,M)\in {\cal S}}M$.
That is, each simple ${\cal A}_{\alpha}(G,{\cal S})$ 
occurs in ${\cal M}$
and its multiplicity space is an irreducible $V^G$-module. 
Moreover,
the different multiplicity spaces are 
inequivalent $V^G$-modules.
It follows {from} this result that
for any $g\in G$ every irreducible $g$-twisted module
is a completely reducible $V^G$-module.
These results are already shown in the case
$g=1$ in \cite{DY} and in the case where 
$g$ is in the center of $G$ in \cite{Y}.

Let ${\cal T}_g$ be the set of all inequivalent irreducible
$g$-twisted $V$-modules for $g\in G$ and set 
${\cal T}=\{(g,M)\ |\ g\in G, M\in {\cal T}_g\}$.
There is a natural right $G$-action on ${\cal T}$. Namely, 
for an irreducible $g$-twisted 
$V$-module $(M, Y_M)$ and $a\in G$, we define  
\[(M, Y_M)\cdot a=(M\cdot a, Y_{M\cdot a}).\]
Here $M\cdot a=M$ as a vector space and $Y_{M\cdot a}(u,z)$ 
is defined by 
\[Y_{M\cdot a}(u,z)=Y_M(au,z)\mbox{ for }u\in V.\]
Note that $M\cdot a$ is an 
irreducible $a^{-1}ga$-twisted $V$-module.
We set $(g,M)\cdot a=(a^{-1}ga, M\cdot a)$ 
for all $(g,M)\in {\cal T}$.

A subset ${\cal S}\subset {\cal T}$
is called {\it stable} if for any $(g,M)\in {\cal S}$ 
and $a\in G$
there exists $(a^{-1}ga, N)\in {\cal S}$ such that 
$M\circ a\simeq N$ as $a^{-1}ga$-twisted $V$-modules.

We assume that ${\cal S}$ is a finite $G$-stable subset 
of ${\cal T}$
until the end of this section.
Let $(g,M)\in {\cal S}$ and $a\in G$.
Then there exists $(aga^{-1},N)\in {\cal S}$
such that $N\cdot a\simeq M$ as $g$-twisted $V$-modules.
That is, there is an isomorphism 
$\phi(a,(g,M)) : M\rightarrow N$
of vector spaces such that 
\[\phi(a,(g,M))Y_{(g,M)}(u,z)
=Y_{(aga^{-1},N)}(au,z)\phi(a,(g,M))\]
for all $u\in V$. By the simplicity of $M$, there exists 
$\alpha_{(g,M)}(a,b)\in {\mathbb C}$
such that 
\[\phi(a,(bgb^{-1},M\cdot b^{-1}))\phi(b,(g,M))
=\alpha_{(g,M)}(a,b)\phi(ab,(g,M)).\]
Moreover, for $a,b,c\in G$ and $(g,M)\in {\cal S}$ we have 
\[\alpha_{(cgc^{-1},M\cdot c^{-1})}(a,b)
\alpha_{(g,M)}(ab, c)=
\alpha_{(g,M)}(ab,c)\alpha_{(g,M)}(b, c).\]
Define a vector space 
${\mathbb C}{\cal S}=\bigoplus_{(g,M)\in {\cal S}}
{\mathbb C}e(g,M)$ with a basis $e(g,M)$ for 
$(g,M)\in {\cal S}$.
The space ${\mathbb C}{\cal S}$ is an associative 
algebra under 
the product $e(g,M)e(h,N)=\delta_{(g,M),(h,N)}(h,N)$.
Let ${\cal U}({\mathbb C}{\cal S})=
\{\sum_{(g,M)\in {\cal S}}\lambda_{(g,M)}e(g,M)\ |
\ \lambda_{(g,M)}\in {\mathbb C}^{\times }\}$
be the set of unit elements on ${\mathbb C}{\cal S}$
where ${\mathbb C}^{\times}$ is the multiplicative group
of ${\mathbb C}$. 
${\cal U}({\mathbb C}{\cal S})$ is a multiplicative right
$G$-module by the action 
$(\sum_{(g,M)\in {\cal S}}\lambda_{(g,M)}e(g,M))\cdot a
=\sum_{(g,M)\in {\cal S}}\lambda_{(g,M)}e(a^{-1}ga,M\cdot a)$ 
for $a\in G$.  
Set $\alpha(a,b)=\sum_{(g,M)\in {\cal S}}
\alpha_{(g,M)}(a,b)e(g,M)$. Then 
\[(\alpha(a,b)\cdot c)\alpha(ab,c)=\alpha(a,bc)\alpha(b,c)\]
holds for all $a,b,c\in G$. 
So $\alpha : G\times G\rightarrow 
{\cal U}({\mathbb C}{\cal S})$ is a 2-cocycle.

Define the vector space
${\cal A}_{\alpha}(G,{\cal S})={\mathbb C}[G]
\otimes {\mathbb C}{\cal S}$ with a basis $a\otimes e(g,M)$ for
$a\in G$ and $(g,M)\in {\cal S}$ and a multiplication on it:
\[a\otimes e(g,M)\cdot b\otimes e(h,N)=
\alpha_{(h,N)}(a,b)ab\otimes e((g,M)\cdot b)e(h,N).\]
Then ${\cal A}_{\alpha}(G,{\cal S})$ is an associative algebra
with the identity element 
$\sum_{(g,M)\in {\cal S}}1\otimes e(g,M)$.

We define an action of ${\cal A}_{\alpha}(G,{\cal S})$ on 
${\cal M}=\bigoplus_{(g,M)\in {\cal S}}M$ as follows:
For $(g,M),(h,N)\in {\cal S},\ w\in N,\ a\in G$ we set 
\begin{eqnarray*}
a\otimes e(g,M)\cdot w & = & \delta_{(g,M),(h,N)}\phi(a,(g,M))w.
\end{eqnarray*}
Note that the actions of ${\cal A}_{\alpha}(G,{\cal S})$
and $V^G$ on ${\cal M}$ commute with each other.

For each $(g,M)\in {\cal S}$ set $G_{(g,M)}
=\{a\in C_G(g)\ |\ M\cdot a\simeq M
\mbox{ as $g$-twisted $V$-modules}\}$. 
Let ${\cal O}_{(g,M)}$
be the orbit of $(g,M)$ under the action of $G$ 
and let $G=\cup_{j=1}^{k}G_{(g,M)}g_j$ 
be a right coset decomposition with $g_1=1$.
Then ${\cal O}_{(g,M)}=\{(g,M)\cdot g_j\ |\ j=1,\ldots,k\}$ and 
${\cal O}_{(g,M)\cdot g_j}=g_j^{-1}G_{(g,M)}g_j$. We define
several subspaces of ${\cal A}_{\alpha}(G,{\cal S})$ by:
\[
\begin{array}{rcl}
{\displaystyle S(g,M)} & = & 
{\displaystyle 
\mbox{Span}\{a\otimes e(g,M)\ |\ a\in G_{(g,M)}\}},\\
{\displaystyle D(g,M)} & = & 
{\displaystyle \mbox{Span}\{a\otimes e(g,M)\ |\ a\in G\}} 
\mbox{  and} \\
{\displaystyle D({\cal O}_{(g,M)}}) & = & 
{\displaystyle \mbox{Span}\{a\otimes 
e(g,M)\cdot g_j\ |\ j=1,\ldots,k, a\in G\}}.
\end{array}
\]
Decompose ${\cal S}$ into a disjoint union of orbits 
${\cal S}=\cup_{j\in J}{\cal O}_j$. Let $(g_j, M^j)$ be a 
representative 
elements of ${\cal O}_j$. Then 
${\cal O}_j=\{(g_j,M^j)\cdot a\ |\ a\in G\}$
and ${\cal A}_{\alpha}(G, {\cal S})
=\bigoplus_{j\in J}D({\cal O}_{(g_j,M^j)})$.
We recall the following properties of 
${\cal A}_{\alpha}(G,{\cal S})$.

\begin{lemma}{\rm (\cite{DY}, Lemma 3.4)}
Let $(g,M)\in {\cal S}$ and $G=\cup_{j=1}^kG_{(g,M)}g_j$. Then
\begin{enumerate}
\item $S(g,M)$ is a subalgebra of ${\cal A}_{\alpha}(G,{\cal S})$
isomorphic to ${\mathbb C}^{\alpha_{(g,M)}}[G_{(g,M)}]$ where 
${\mathbb C}^{\alpha_{(g,M)}}[G_{(g,M)}]$ is the twisted 
group algebra with 
$2$-cocycle $\alpha_{(g,M)}$.
\item $D({\cal O}_{(g,M)})=\oplus_{j=1}^{k}D((g,M)\cdot g_j)$ 
is a direct sum of left ideals.
\item Each $D({\cal O}_{(g,M)})$ is a two sided ideal of 
${\cal A}_{\alpha}(G, {\cal S})$ and 
${\cal A}_{\alpha}(G, {\cal S})
=\oplus_{j\in J}D({\cal O}_{(g_j,M^j)})$.
Moreover, $D({\cal O}_{(g,M)})$ has the identity element 
$\sum_{(h,N)\in {\cal O}_{(g,M)}}1\otimes e(h,N)$.
\end{enumerate}
\end{lemma}

\begin{lemma}{\rm (\cite {DY}, Theorem 3.6)}\label{Lemma : rep} 
\begin{enumerate}
\item $D({\cal O}_{(g,M)})$ is semisimple for all 
$(g,M)\in {\cal S}$
and the simple $D({\cal O}_{(g,M)})$-modules are precisely equal to 
$\mbox{\rm Ind}^{D(g,M)}_{S(g,M)}W=D(g,M)\otimes_{S(g,M)}W$ 
where $W$ ranges over the simple 
${\mathbb C}^{\alpha_{(g,M)}}[G_{(g,M)}]$-modules.
\item ${\cal A}_{\alpha}(G, {\cal S})$ 
is semisimple and simple ${\cal A}_{\alpha}(G, {\cal S})$-modules 
are precisely
$\mbox{\rm Ind}^{D(g_j,M^j)}_{S(g_j,M^j)}W$ 
where $W$ ranges over the simple 
${\mathbb C}^{\alpha_{(g_j,M^j)}}[G_{(g_j,M^j)}]$-modules 
and $j\in J$.
\end{enumerate}
\end{lemma}

Let $(g,M)\in {\cal T}$. 
Following \cite{Z} we define weight zero operator $o_M(u)$ by 
$u_{\mwt{u}-1}$ on $M$ for homogeneous $u\in V$
and extend $o_M(u)$ to all $u$ by linearity.
For any nonnegative rational number $n\in {\mathbb Z}/T$, 
define a map $\sigma_n$ {from} $A_{G,n}(V)$ 
to ${\displaystyle \bigoplus_{(g,M)\in {\cal S}}
\bigoplus_{\stackrel{\scriptstyle m\in {\mathbb Z}/T}{0\leq m\leq n}}
\mbox{End}M(m)}$ by $\sigma_n(u)=\sum_{(g,M)\in {\cal S}}
\sum_{\stackrel{\scriptstyle m\in {\mathbb Z}/T}{0\leq m\leq n}}
o(u)_M$. For $a\in G$ and $f\in \mbox{End}(M)$,
define the action of $a$ by 
$a\cdot f=\phi(a,(g,M))f\phi(a,(g,M))^{-1}
\in \mbox{End}(M\cdot a^{-1})$.
This defines a left action of $G$ on 
${\displaystyle \bigoplus_{(g,M)\in {\cal S}}
\bigoplus_{\stackrel{\scriptstyle m\in 
{\mathbb Z}/T}{0\leq m\leq n}}
\mbox{End}M(m)}$.

We prepare the following results in order to show the main result
in this section.
For any $(g,M)\in {\cal S}$, we always arrange the grading 
on $M=\oplus_{0\leq j\in {\mathbb Z}/T}M(j)$ so that
$M(0)\neq 0$ if $M\neq 0$ using a grading shift.

\begin{lemma}\label{Lemma : doukei}
Let $(g,M),(h,N)\in {\cal S}$. Let $m\in {\mathbb Z}/|g|,
n\in {\mathbb Z}/|h|$ such that $m\leq n$.
If $M(m)\neq 0$ and 
$M(m)\simeq N(n)$ as $A_{G,n}(V)$-modules, then
$(g,M)=(h,N)$ and $m=n$.
\end{lemma}
\pf
Since $M(m)\neq 0$ and 
$M(m)\simeq N(n)$ as $A_{G,n}(V)$-modules, 
we have 
$g=h$ {from} Theorem \ref{Theorem : zhu}.
So $M(m)\simeq N(n)$ as $A_{g,n}(V)$-modules.
We have $m=n$ and $M=N$ by Theorem 4.3 in \cite{DLM2}. $\Box$

\begin{lemma}\label{Lemma : density}
The map $\sigma_n$ is a $G$-module epimorphism. In particular,
\[{\displaystyle \sigma_n(A_{G,n}(V)^G)
=(\bigoplus_{(g,M)\in {\cal S}}
\bigoplus_{\stackrel{\scriptstyle m\in 
{\mathbb Z}/T}{0\leq m\leq n}}
\mbox{\rm End}M(m))^G}.\]
\end{lemma}
\pf The proof is similar to that of in Lemma 6.12 in \cite{DY}
because of Lemma \ref{Lemma : doukei}. $\Box$ \\

For $(g,M)\in {\cal S}$ 
let $\Lambda_{G_{(g,M)}, \alpha_{(g,M)}}$ be the set of 
all irreducible characters $\lambda$ of 
${\mathbb C}^{\alpha_{(g,M)}}[G_{(g,M)}]$. We denote the 
corresponding
simple module by $W(g,M)_{\lambda}$. Note that $M$ is a
semisimple ${\mathbb C}^{\alpha_{(g,M)}}[G_{(g,M)}]$-module.
Let $M^{\lambda}$ be the sum of simple 
${\mathbb C}^{\alpha_{(g,M)}}[G_{(g,M)}]$-module of $M$ 
isomorphic to $W(g,M)_{\lambda}$. Then 
$M=
\oplus_{\lambda\in \Lambda_{G_(g,M),\alpha_{(g,M)}}}M^{\lambda}$.
Moreover $M^{\lambda}=W(g,M)_{\lambda}\otimes M_{\lambda}$
where $M_{\lambda}=
\mbox{Hom}_{{\mathbb C}^{\alpha_{(g,M)}}[G_{(g,M)}]}
(W(g,M)_{\lambda},M)$ is the multiplicity of $W(g,M)_{\lambda}$ 
in $M$.
We can realize $M_{\lambda}$ as a subspace of $M$
in the following way: Let $w\in W(g,M)_{\lambda}$ be a 
fixed nonzero
vector. Then we can identify 
$\mbox{Hom}_{{\mathbb C}^{\alpha_{(g,M)}}[G_{(g,M)}]}
(W(g,M)_{\lambda},M)$  with the subspace 
\[\{f(w)\ |\ f\in 
\mbox{Hom}_{{\mathbb C}^{\alpha_{(g,M)}}[G_{(g,M)}]}
(W(g,M)_{\lambda},M)\},\]
of $M^{\lambda}$. Note that the actions of 
${\mathbb C}^{\alpha_{(g,M)}}[G_{(g,M)}]$
and $V^{G_{(g,M)}}$ on $M$ commute with each other. So
$M^{\lambda}$ and $M_{\lambda}$ are ordinary 
$V^{G_{(g,M)}}$-modules. 
Furthermore, 
$M^{\lambda}$ and $M_{\lambda}$ are ordinary $V^{G}$-modules. 
It is shown by Theorem 3.9 in \cite{Y} that 
$M_{\lambda}$ is a 
nonzero irreducible $V^{G_{(g,M)}}$-module for any
$\lambda\in \Lambda_{G_{(g,M)}, \alpha_{(g,M)}}$.

Let ${\cal S}=\cup_{j\in J}{\cal O}_j$ be 
an orbit decomposition and fix 
$(g_j,M^j)\in {\cal O}_j$ for each $j\in J$. 
For convenience, we set 
$G_j=G_{(g_j,M^j)}, 
\Lambda_j=\Lambda_{(g_j,M^j),\alpha_{(g_j,M^j)}}$ and
$W_{j, \lambda}=W(g_j,M^j)_{\lambda}$ for $j\in J$ and
$\lambda\in \Lambda_j$.
We have a decomposition
\[M^j=\bigoplus_{\lambda\in \Lambda_j}
W_{j,\lambda}\otimes M^j_{\lambda},\]
as a
${\mathbb C}^{\alpha_{(g_j,M^j)}}[G_{j}]
\otimes V^{G_{j}}$-module.
%where $M^j_{\lambda}$ is an irreducible $V^{G_j}$-module. 
We also have
\[{\cal M}=\bigoplus_{j\in J, \lambda\in \Lambda_j}
\mbox{Ind}^{D(g_j,M^j)}_{S(g_j,M^j)}
W_{j,\lambda}\otimes M^{j}_{\lambda}.\]
as a ${\cal A}_{\alpha}(G, {\cal S})\otimes V^G$-module using 
the same arguments as those in Proposition 6.5 in \cite{DY}.
 
For $j\in J$ and
$\lambda\in \Lambda_j$ we set $W^{j}_{\lambda}=
\mbox{Ind}^{D(g_j,M^j)}_{S(g_j,M^j)}W_{j, \lambda}$.
Then $W^{j}_{\lambda}$ forms a complete list of simple 
${\cal A}_{\alpha}(G, {\cal S})$-modules {from} 
Lemma \ref{Lemma : rep}.
Now we have the main result in this section.

\begin{theorem}\label{theorem : main}
As a ${\cal A}_{\alpha}(G,{\cal S})\otimes V^G$-module,
\[{\cal M}=\bigoplus_{j\in J,\lambda\in\Lambda_j}
W^{j}_{\lambda}\otimes M^j_{\lambda}.\]
Moreover
\begin{enumerate}
\item Each $M^{j}_{\lambda}$ is a nonzero 
irreducible $V^G$-module.
\item $M^{j_1}_{\lambda_1}$ and $M^{j_2}_{\lambda_2}$ 
are isomorphic 
$V^G$-module if and only if $j_1=j_2$ and $\lambda_1=\lambda_2$.
\end{enumerate}
In particular, all irreducible $g$-twisted $V$-modules are 
completely irreducible $V^G$-modules. 
\end{theorem}
\pf The proof is similar to that of Theorem 6.13 in \cite{DY}
because of Lemma \ref{Lemma : density}.$\Box$

\section{Appendix}
\setcounter{equation}{0}
In this appendix we prove the matrix $A_3$ in (\ref{eqn : det2})
is non-singular. Let $a, b, t$ positive integers and 
$x_1,x_2,\ldots,x_t$ indeterminants.
Set a $bt\times bt$-matrix
\begin{eqnarray*}
& & A = \\
& & \left[
\begin{array}{ccccccccc}
{x_1\choose a}& {x_1\choose a+1}&\cdots & {x_1\choose a+b-1}&
{x_2\choose a}& \cdots & {x_2\choose a+b-1}&
\cdots & {x_t\choose a+b-1}\\
{x_1\choose a+1}& {x_1\choose a+2}&\cdots & {x_1\choose a+b}&
{x_2\choose a+1}& \cdots & {x_2\choose a+b}&
\cdots & {x_t\choose a+b}\\
\vdots & \vdots & & \vdots & \vdots & & \vdots & & \vdots\\
{x_1\choose a+bt-1}& {x_1\choose a+bt}&\cdots & 
{x_1\choose a+bt+b-2}&
{x_2\choose a+bt-1}& \cdots & {x_2\choose a+bt+b-2}&
\cdots & {x_t\choose a+bt+b-2}
\end{array}\right].
\end{eqnarray*}
It is proved {from} the following result 
that the matrix $A_3$ is non-singular.

\begin{proposition}
Let 
\begin{eqnarray*}
H(x_1,\ldots,x_t) & = & 
\prod_{i=1}^{t}
\left(\prod_{j=1}^{b-1}(x_i+j)^{b-j}
\prod_{j=0}^{a-1}(x_i-j)^b
\prod_{j=1}^{b-1}(x_i-a+1-j)^{b-j}\right)\\
& & \times 
\prod_{1\leq i<j\leq t}
\prod_{k=-b+1}^{b-1}(x_i-x_j+k)^{b-|k|}.
\end{eqnarray*}
Then
\begin{eqnarray*}
\det A & = & (-1)^{bt(bt-1)/2}
\frac{H(x_1,\ldots,x_t)}
{H(a+tb-1,a+(t-1)b-1,\ldots,a+b-1)}.
\end{eqnarray*}
\end{proposition}
\pf
Since ${x\choose j}=x(x-1)\cdots (x-j+1)/j!$, 
$\prod_{j=0}^{a-1}(x_i-j)^b 
\prod_{j=1}^{b-1}(x_i-a+1-j)^{b-j}$ is a factor
of $\det A$ for any $1\leq i\leq t$.
For any $1\leq i\leq t$ and $0\leq k\leq b-1$, set
a $bt\times b$-matrix $B(i,k)$ by  
\begin{eqnarray*}
B(i,k) & = & 
\left[
\begin{array}{ccccccc}
{x_i\choose a} & {x_i+1\choose a+1} 
& \cdots & {x_i+k\choose a+k} 
& {x_i+k\choose a+k+1} & \cdots & {x_i+k\choose a+b-1}\\
{x_i\choose a+1} & {x_i+1\choose a+2} 
& \cdots & {x_i+k\choose a+k+1} 
& {x_i+k\choose a+k+2} & \cdots & {x_i+k\choose a+b}\\
\vdots & \vdots & & \vdots & \vdots & & \vdots \\
{x_i\choose a+bt-1} & {x_i+1\choose a+bt}& \cdots & 
{x_i+k\choose a+bt-1+k} 
& {x_i+k\choose a+bt+k} & \cdots & {x_i+k\choose a+bt+b-2}\\
\end{array}\right].
\end{eqnarray*}
For any $1\leq k\leq b-1$,
$\det A$ is equal to that of the $bt\times bt$-matrix
$\left[B(1,k)\ B(2,k)\cdots B(t,k)\right]$
because ${x\choose j}+{x\choose j+1}={x+1\choose j+1}$.
So $(x_i+k)^{b-k}$ is a factor
of $\det A$ for any $1\leq i\leq t$ and $1\leq k\leq b-1$.
Fix any $1\leq i<j\leq t$ and $0\leq k\leq b-1$.
Then $\det A$ is equal to that of the $bt\times bt$-matrix
\[
\left[B(1,0)\ \cdots\ B(i,0)\ \cdots\ B(j,k)\ \cdots
\ B(t,0)\right]
\]
by the same reason as above. 
Comparing the $(t(i-1)+p)$-th column and 
the $(t(j-1)+p)$-th column
for all $k+1\leq p\leq b$, we have 
$(x_i-x_j-k)^{b-k}$ is a factor of $\det A$.
We also have that $(x_i-x_j+k)^{b-k}$ is a factor of 
$\det A$ by applying the same argument to the matrix
\[
\left[B(1,0)\ \cdots\ B(i,k)\ \cdots\ B(j,0)\ \cdots
\ B(t,0)\right].
\]
So there is $\alpha(a,b)\in {\mathbb C}[x_1,\ldots,x_t]$
such that  
\begin{eqnarray}
\det A & = & \alpha(a,b) H(x_1,\ldots,x_n)\label{eqn : han}
\end{eqnarray}
We have $\alpha(a,b)\in {\mathbb C}$
since the degrees of both sides in formula 
(\ref{eqn : han}) are equal to $(a-1)bt+b^2t(t+1)/2$. 
Substituting $(a+tb-1,a+(t-1)b-1,\ldots,a+b-1)$
for $(x_1,x_2,\ldots,x_t)$ in $A$, we have an anti-diagonal matrix
with all anti-diagonal elements $1$. Hence 
$\alpha(a,b)=(-1)^{bt(bt-1)/2}/H(a+tb-1,a+(t-1)b-1,\ldots,a+b-1)$.
$\Box$

\bigskip
\noindent{\bf Acknowledgments} 

\noindent We would like to thank Soichi Okada 
for information of a method to compute the determinant of 
the matrix in Appendix.

\end{document}